\documentclass[runningheads, 12pt]{llncs}
\usepackage{amssymb}
\usepackage [english]{babel}
\usepackage[utf8]{inputenc}
\usepackage[T1]{fontenc}
\usepackage{caption}
\usepackage{amssymb}
\usepackage{amsmath}
\usepackage{amsfonts}
\usepackage{booktabs}
\usepackage{graphicx}
\usepackage{url}
\usepackage[usenames]{color}
\usepackage{algorithm}
\usepackage{algpseudocode}
\usepackage{parskip}
\usepackage{makecell}
\usepackage{xfrac}
\usepackage{xifthen}
  \usepackage{geometry}
\newcommand{\Nedelec}{N\'{e}d\'{e}lec }

\newcommand{\keywords}[1]{{\scriptsize Keywords: #1}}
\newcommand{\sep}{$\cdot$ }

\begin{document}

\urldef{\mailsa}\path|stefan.takacs@ricam.oeaw.ac.at|

\title{Robust multigrid methods for isogeometric discretizations of the Stokes equations}
\author{Stefan Takacs$^1$}
\institute{ $^1$ Johann Radon Institute for Computational and Applied Mathematics (RICAM),\\
Austrian Academy of Sciences, Altenberger Str. 69, 4040 Linz, Austria\\
\mailsa
 }

\noindent
\maketitle

\begin{abstract}
In recent publications, the author and his coworkers have proposed a multigrid
method for solving linear systems arizing from the discretization of partial
differential equations in isogeometric analysis and have proven that the
convergence rates are robust in both the grid size and the polynomial degree.
So, far the method has only been discussed for the Poisson problem. 
In the present paper, we want to face the question if it is possible to extend the
method to the Stokes equations.
\end{abstract}

\keywords{
Isogeometric analysis \sep Multigrid methods \sep Stokes problem
}


\pagestyle{myheadings}
\thispagestyle{plain}
\markboth{}{S. Takacs, Multigrid for IgA-Stokes}

\section{Introduction}

Isogeometric analysis (IgA) was introduced by Tom Hughes et al. in the 2005, cf.~\cite{Hughes:2005}, aiming
to improve the connection between computer aided design (CAD) and finite element (FEM) simulation.
In IgA, as in CAD software, B-splines and non-uniform rational B-splines (NURBS) are used for representing
both the geometrical objects of interest and the solution of the partial differential equation (PDE) to be
solved.

In IgA, mostly B-splines or NURBS of maximum smoothness are used, i.e., having a spline degree of $p$,
the functions are $p-1$ times continuously differentiable. Using such a function space, one obtains on the one
hand the approximation power of high order functions, while on the other hand, unlike in standard
high-order FEM, one does not suffer from a growth of the number of degrees of freedom.

From the computational point of view, the treatment of the linear systems arizing from the discretization
with high spline degrees is still challenging as the condition number both of mass and stiffness matrices
grows exponentially with the spline degree. In the early IgA literature, often 
finite element solvers have been transferred to IgA only with minimal
adaptations. Numerical experiments indicate that such approaches result in methods
that work well for small spline degrees, but their performance deteriorates as the
degree is increased, often dramatically.  In~\cite{HTZ:2016,HT:2016}, the author and his coworkers have
proposed multigrid methods which are provable robust in the polynomial degree and the grid
size. Numerical experiments indicated that the proposed approach of \emph{subspace corrected mass
smoothers} seems to pay off (compared to multigrid methods with a standard Gauss-Seidel smoother)
for polynomial degrees of four or five.

In the present paper, we want to discuss the extension of the subspace corrected mass smoothers
beyond the case of the Poisson problem to the Stokes flow problem. Unlike for the Poisson problem,
for the Stokes problem already the setup of a stable isogeometric discretization is non standard. As there have already
been results in the literature, we refer to paper~\cite{BuffaEtAl}, which serves as a basis of the present
paper.  Alternative approaches
can be found in~\cite{Evans:2013,Bressan:2013,Buffa:2011} and others.
After introducing discretizations, we discuss the setup of the preconditioner.

For the Poisson problem, the multigrid solver has been applied directly and as a preconditioner
for the conjugate gradient method. For the case of a non-trivial geometry transformation,
in \cite{HT:2016} a conjugate gradient method, preconditioned with the multigrid method for the
parameter domain, has been used. It has been shown that in this case the resulting method is robust both the grid size and the polynomial
degree, but not in the geometry transformation.

There are a few approaches how to carry this over to the Stokes equations. The first possibility
is to apply the multigrid method directly to the problem of interest (all-at-once multigrid method), cf.~\cite{Vanka:1986} for
a particularly popular method in standard FEM or~\cite{Larin:Reusken:2008} for a survey.
As the results for the Poisson problem have indicated that a direct application of the multigrid method
in the presence of a non-trivial geometry transformation is not optimal, we do not concentrate on that case.

So, we consider a Krylov space method with an appropriate preconditioner, living on the parameter domain. In principle,
this could be the Stokes problem on the physical domain, but such a choice (an indefinite preconditioner for an indefinite
problem) typically requires the use of a GMRES method, whose convergence cannot be easily proven by considering
the spectrum of the preconditioned system, cf.~\cite{Greenbaum:1996}.
So, we consider elliptic preconditioners, particularly block-diagonal preconditioners.
As the Stokes equations are well-posed in the Sobolev space $H^1$ (velocity) and the Lebesgue space $L_2$ (pressure), we just
setup preconditioners for those spaces (operator preconditioning).

As we observe that the subspace corrected mass smoothers suffer significantly from the geometry transformation,
we propose a variant (by incorporating an approximation to the geometry transformation) which leads for
several experiments to a significant speedup.

This paper is organized as follows. We will introduce the
particular model problem in section~\ref{sec:1} and discuss three kinds of discretizations for the mixed
system in section~\ref{sec:2}. As a next step, in section~\ref{sec:3}, we propose a preconditioner.
Finally, in section~\ref{sec:4}, we give the results of the numerical examples and draw some conclusions.

\section{Model problem}\label{sec:1}

Let $\Omega\subseteq \mathbb{R}^2$ be a simply connected domain with Lipschitz
boundary~$\partial\Omega$ and assume a force filed $f$ to be given on~$\Omega$ and boundary data
to be given on~$\partial\Omega$. The \emph{Stokes flow model problem} reads as follows.
Find the velocity field $u$ and the pressure distribution $p$ such that
\begin{equation}\label{eq:stokes}
        -\Delta u + \nabla p = f \qquad\mbox{and} \qquad \qquad \nabla \cdot u = 0
\end{equation}
hold on $\Omega$ and Dirichlet boundary conditions hold on $\partial\Omega$.
After homogenization, we obtain a mixed variational form, which reads as follows.
Find~$u\in V:=H^1_0(\Omega)$
and~$p\in Q:=L_2(\Omega)$ such that
\[
      \underbrace{( \nabla u, \nabla v)}_{\displaystyle a(u,v)}
    + \underbrace{( \nabla\cdot v, p)}_{\displaystyle b(v,p)}
    = (f,v) \quad\forall v\in V,
  \qquad
      \underbrace{( \nabla\cdot u, q)}_{\displaystyle b(u,q)}
    = 0 \quad\forall q\in Q.
\]
Here, and in what follows $L_2(\Omega)$, $H^1(\Omega)$ and $H^1_0(\Omega)$ are the
standard Lebesgue and Sobolev spaces, and $(\cdot,\cdot)$ is the standard norm on $L_2(\Omega)$.

Existence and uniqueness of the solution and its dependence of the data follows from Brezzi's theorem~\cite{Brezzi:1974},
which requires besides boundedness and $H^1$-coercivity of~$a$ the
inf-sup stability
\begin{align*}
    \inf_{q \in L_2(\Omega)}
    \sup_{v \in H^1(\Omega)}
       \frac{  ( \nabla\cdot v, q) }{ \| v\|_{H^1(\Omega)} \|q\|_{L_2(\Omega)}   }
    \ge C,
\end{align*}
which is known to be satisfied for the Stokes problem, cf.~\cite{Brezzi:Fortin:1991}.

\section{Discretization}\label{sec:2}

The discretization is done using a standard Galerkin approach, i.e., we replace the spaces $V$ and
$Q$ by finite-dimensional subspaces $V_h$ and $Q_h$.
As for the continuous problem, existence and uniqueness of the solution can be shown by Brezzi's theorem.
Boundedness and  $H^1$-coercivity of $a$ follow directly from the continuous problem,
but the inf-sup stability for the discrete problem does not. Therefore, we have to guarantee that the
\emph{discrete inf-sup condition}
\begin{align*}
    \inf_{q_h \in Q_h}
    \sup_{v_h \in V_h}
        \frac{  ( \nabla\cdot v_h, q_h) }{ \| v_h\|_{H^1(\Omega)} \|q_h\|_{L_2(\Omega)}   }
    \ge C
\end{align*}
is satisfied, which is actually a condition on the discretization. In the subsection~\ref{subsec:2:2},
we will discuss discretizations satisfying this condition.

Assuming a particular discretization and a basis for the chosen space, one ends
up with a linear system to be solved: For given $\underline{f}_h$, find $\underline{x}_h$ such that
\begin{equation}\label{eq:linsys}
    A_h\, \underline{x}_h = \underline{f}_h,
        \qquad\mbox{where}\qquad
        A_h = \left(    \begin{array}{cc}   K_h & D_h^T \\ D_h  &0     \end{array}       \right)
        \qquad \mbox{and} \qquad
        \underline{x}_h = \left(    \begin{array}{c}  \underline{u}_h \\ \underline{p}_h       \end{array}       \right)
\end{equation}
and $K_h$ is a standard stiffness matrix and $D_h$ is a matrix representing the divergence.

\subsection{Discretization in isogeometric analysis}\label{subsec:2:1} 

Let $S_{p,h}^q$ be the
space of all $q$ times continuously differentiable functions on $(0,1)$, which are piecewise polynomials
of degree $p$ on a (uniform) grid of size $h=1/n$.
As a basis for $S_{p,h}^q$ we choose the classical basis of B-splines, see, e.g., \cite{DeBoor:1972}.

For computational domains $\Omega\subset \mathbb{R}^2$, we first define the
spline spaces for the \emph{parameter domain} $\hat{\Omega} = (0,1)^2$.
On the parameter domain, we introduce the space of tensor-product splines, which reads as follows:
\[S_{p_1,p_2,h}^{q_1,q_2}:= S_{p_1,h}^{q_1}\otimes
S_{p_2,h}^{q_2},\] where $A\otimes B$ denotes the linear span of all functions
$(x,y)\mapsto u(x) v(y)$, where $u\in A $ and $v\in B$.
Note that the restriction to two dimensions and to a uniform grid is only for ease
of notation. The extension to three and more dimensions or to non-uniform grids is
completely straight-forward.
Assuming that \emph{physical domain} $\Omega$ is the image of a B-spline or NURBS
mapping
\[
G: \hat{\Omega}=(0,1)^2 \rightarrow \Omega,
\]
we define the spline spaces on the physical domain typically using a classical pull-back principle.
More complicated domains are represented patch-wise, where for each patch a separate geometry
transformation $G$ exists. For simplicity, we do not discuss that in the present paper.

\subsection{Stable discretizations for the Stokes problem}\label{subsec:2:2} 

As mentioned above, it is required to set up the discretization such that the discrete inf-sup condition
holds. We discuss this first for the parameter domain. Here, we follow the outline of the
paper~\cite{BuffaEtAl}, where three spline space configurations
have been proposed, which are variants of known stable spaces from standard finite elements:
Taylor-Hood like splines~$\hat{X}_h^{(\textnormal{TH})}$, \Nedelec like splines~$\hat{X}_h^{(\textnormal{NE})}$
and Raviart-Thomas like splines~$\hat{X}_h^{(\textnormal{RT})}$.
All of them utilize the same grid for both the velocity and the pressure, which makes the implementation
significantly easier compared to approaches that are based on setting up two different grids (like IgA-variants
of the macro elements, cf.~\cite{Bressan:2013}). All of these discretizations follow the spirit of IgA, allowing to freely choose the underlying
polynomial degree~$p$. For all of them, the smoothness is in the order of the polynomial degree, which preserves
the feature that the number of degrees of freedom is basically not increased when the polynomial degree is increased.

For the case of two dimensions, the spaces are given by
\begin{equation*}
    \begin{array}{lll}
   \hat{X}_h^{(\textnormal{TH})} := \hat{V}_h^{(\textnormal{TH})}  \times \hat{Q}_h, 
  & \qquad\hat{V}_h^{(\textnormal{TH})} := S_{p+1,p+1}^{p-1,p-1}\times S_{p+1,p+1}^{p-1,p-1},
  &   \\
   \hat{X}_h^{(\textnormal{NE})} := \hat{V}_h^{(\textnormal{NE})}  \times \hat{Q}_h,
  & \qquad\hat{V}_h^{(\textnormal{NE})} := S_{p+1,p+1}^{p,p-1}\times S_{p+1,p+1}^{p-1,p}, 
  &  \\
   \hat{X}_h^{(\textnormal{RT})} := \hat{V}_h^{(\textnormal{RT})}  \times \hat{Q}_h,
  & \qquad\hat{V}_h^{(\textnormal{RT})} := S_{p+1,p}^{p,p-1}\times S_{p,p+1}^{p-1,p}, 
  & \qquad\hat{Q}_h := S_{p,p}^{p-1,p-1},\\
  \end{array}
\end{equation*}
where $A\times B:=\{(a,b):a\in A,\;b\in B\}.$
Observe that these spline spaces are nested, i.e., we have
    \[    \hat{V}_h^{(\textnormal{RT})} \subset \hat{V}_h^{(\textnormal{NE})} \subset \hat{V}_h^{(\textnormal{TH})}\]
and (for $n>\hspace{-.3em}>p$) a ratio of $9:5:3$ for the number of degrees of freedom.
The extension of these definitions to three and more dimensions is straight-forward, cf.~\cite{BuffaEtAl}.

For all of these settings, the discrete inf-sup condition has been shown in~\cite{BuffaEtAl}. For the
Raviart-Thomas like splines, the discrete inf-sup condition cannot be proven if Dirichlet boundary conditions
are present. As the method still seems to work well in practice, we include also the Raviart-Thomas
discretization in our experiments.

The next step is to introduce the discretization for the physical domain.
As outlined in the beginning of this section, the discretization, once introduced
on the parameter domain, is typically defined on the physical domain just by
\emph{direct composition}:
\[
        V_h^{(X,\textnormal{D})} := \{ v_h \;|\; v_h \circ G \in \hat{V}_h^{(X)} \},
        \qquad X\in\{\textnormal{TH},\textnormal{NE},\textnormal{RT}\}.
\]
For the Stokes problem, as an alternative, the divergence preserving 
\emph{Piola transform} has been proposed:
\[
        V_h^{(X,\textnormal{P})} := \left\{ v_h \;\left|\; \frac{1}{\mbox{det}\,J_G} J_G\,v_h \circ G \in \hat{V}_h^{(X)} \right.\right\},
        \qquad X\in\{\textnormal{TH},\textnormal{NE},\textnormal{RT}\},
\]
where $J_G$ is the Jacobi matrix of $G$.
The pressure distribution, which is a scalar quantity, is always mapped directly, i.e., in all cases
we choose the direct composition
\[
        Q_h := \{ q_h \;|\; q_h \circ G \in \hat{Q}_h \}.
\]
In~\cite{BuffaEtAl}, the inf-sup stability has been shown if the Piola transform is used, for the
Taylor-Hood like splines also if the direct composition is used.
Again, we report also on the numerical results for the cases that are not covered by the
convergence theory (direct composition for the \Nedelec like and the Raviart-Thomas like splines).

\section{Robust multigrid solvers}\label{sec:3}

As outlined in the introduction, the multigrid preconditioner aims to represent the
theoretical block-diagonal preconditioner
\[
        Q_h := \left(    \begin{array}{cc}   K_h &  \\  & \beta^{-1} M_h     \end{array}       \right),
\]
where $K_h$ is the stiffness matrix, $M_h$ is the mass matrix and $\beta>0$ is an
accordingly chosen scaling parameter. As mentioned above
and as discussed in detail in~\cite{HT:2016}, we use as preconditioner
for the problem on the physical domain the corresponding preconditioner, say $\hat{Q}_h$, on the parameter domain.
There, the matrices $M_h$ and $K_h$ are replaced by $\hat{M}_h$ and $\hat{K}_h$, their
counterparts on the parameter domain. Note that the stiffness matrix acts on the velocity
variable, a vector-valued quantity, and that this matrix is block-diagonal on the parameter
domain and, iff the direct composition is used, on the physical domain.
In all cases, $K_h$ and $\hat{K}_h$ are spectrally equivalent.

Instead of an exact inverse of the matrix $\hat{Q}_h$, we only need to realize
an approximation to the application of $\hat{K}_h^{-1}$ and $\hat{M}_h^{-1}$ to any given vector.
The approximation of  $\hat{K}_h^{-1}$ is realized using one multigrid
V-cycle with one pre- and one post-smoothing step of the subspace corrected mass smoother,
as proposed in~\cite{HT:2016}. There, the algorithm was analyzed only for the
case of splines of maximum smoothness, however it can be applied for any spline space
and robustness in the polynomial degree can be guaranteed by a slight extension of the presented
theory as long as the smoothness is in the order of the polynomial degree.
As in the previous publications~\cite{HTZ:2016,HT:2016}, the grid hierarchy 
is set up for a fixed polynomial degree and a fixed smoothness by just uniformly refining the grid.
Using this approach, one obtains nested spaces, so the setup of the coarse-grid correction is trivial.

One of the key observations which was leading to the results in~\cite{HTZ:2016,HT:2016} was that the
spectral equivalence of the mass matrix and its diagonal deteriorates if $p$ is increased. This
has also to be taken into account when constructing the preconditioner for the pressure variable.
Analogously to the smoother, we realize
the application of $\hat{M}_h^{-1}$ exactly, based on the tensor-product structure
of the mass matrix.

The preconditioner is symmetric and positive definite and can therefore be applied in
the framework of a MINRES iteration.

\section{Numerical results}\label{sec:4}

The numerical experiments have been performed using the C++ library G+SMO, see~\cite{gismoweb}, both for the unit square,
i.e., for a problem without geometry transformation, and for a quarter annulus $\{(x,y)\in\mathbb{R}_+^2:1< x^2+y^2 < 4\}$.
For both problems, the problem has been constructed (with inhomogeneous right-hand-side and inhomogeneous
Dirichlet boundary conditions) such that the exact solution is 
\[ 
    u_h(x,y)= \left(\begin{array}{c}
     \cos(5x+5y)+\sin(5x-5y)\\
    -1-\cos(5x+5y)+\sin(5x-5y)
    \end{array}\right),
\]
and $p_h(x,y) = -(1+x)(1+y)+c$, where $c$ is chosen such that $\int_\Omega p_h\, \mbox{d}x=0$.

In Table~\ref{tab:paramdom}, we see number of MINRES steps required for reducing
the initial error (measured in the $\ell^2$-norm of the solution vector) by a factor of
$10^{-6}$; cases where the memory was not enough are indicated with OoM. The discussion
is done for all proposed discretization schemes. 
The need of the discussion of $p$-robust methods is easily observed
when looking at the results for a standard preconditioner: We display the results if
one multigrid V-cycle with Gauss-Seidel smoother is used for the velocity and one symmetric Gauss-Seidel
sweep is used for the pressure (\emph{GS-MG}). There, the number of iterations increases
drastically if $p$ is increased. As the approach is perfectly robust in the grid size~$h=2^{\ell}$,
we omit the numbers for finer grids.
Compared to that approach, the preconditioner proposed in Section~\ref{sec:3} (\emph{SCMS-MG})
shows results which
are robust both in the grid size and the polynomial degree and which works well for all discretizations.
Although the iteration numbers are smaller than for the GS-MG preconditioner, one has to consider that
the costs of the SCMS-MG preconditioner are significantly higher than those of the GS-MG preconditioner,
so the proposed method only pays of iff higher polynomial degrees (starting from $4$ or $5$)
are considered. We have chosen $\beta=0.05$ and as damping parameter $\sigma$ of the underlying smoother, cf.~\cite{HT:2016}, either
$\sigma^{-1} = 0.04 \, \hat{h}^2$ (for Taylor-Hood and N\'{e}d\'{e}lec) or $0.16 \,\hat{h}^2$ (for Raviart-Thomas),
where $\hat{h}$ is the grid size on the parameter domain. While some of the numbers
might be improved by fine-tuning the parameters, the given
tables for reasonable uniform choices show what one can expect for each of the methods.

In Table~\ref{tab:approx} we see how well the computed solution approximates the exact solution
in terms of the $L^2$-norm. Here, we have used the abovementioned solver, where the
the stopping criterion has been chosen to reach either a relative error of $10^{-10}$ or $100$ iterations.
We present the error between the computed solution and the known exact solution (for the pressure after
projecting into the space of functions with vanishing mean). 
We observe that, for the same choice of the polynomial order $p$ and the same
grid size, the Taylor-Hood discretization yields to the smallest errors,
for the cost of the largest number of degrees of freedom. For the Raviart-Thomas discretization (where the
inf-sup condition cannot be shown for the chosen Dirichlet boundary conditions), we observe that the
error for the velocity converges, while the error of the pressure stagnates at around  $10^{-2}$.
Observe moreover that for $p=5$, the approximation on the coarsest grid was fine enough such that
the approximation error could not be improved by refinement.

\newcommand{\dnc}{{\tiny >1k}}
\newcommand{\oom}{{\tiny OoM}}

\begin{table}[p]
   \centering
    \begin{tabular}{p{.3em}p{2.2em}|p{.1cm}p{2em}p{2em}p{2em}p{2em}|p{.1cm}p{2em}p{2em}p{2em}p{2em}|p{.1cm}p{2em}p{2em}p{2em}p{2em}}
    \toprule
         &           & \multicolumn{5}{c|}{Taylor-Hood} &  \multicolumn{5}{c|}{\Nedelec} & \multicolumn{5}{c}{Raviart-Thomas} \\
    \quad & $\ell\,\diagdown\, p$     &&  2 & 3 & 5 & 8 && 2 & 3 & 5 & 8 && 2 & 3 & 5 & 8    \\
    \midrule
     \multicolumn{17}{l}{MINRES, preconditioned with SCMS-MG} \\
    \midrule
     & 5 &&  55 & 54 & 49 & 46         &&  80 & 74 & 68 & 55          &&  44 & 36 & 35 & 29 \\
     & 6 &&  54 & 58 & 53 & 51         &&  76 & 76 & 70 & 63          &&  44 & 37 & 36 & 32 \\
     & 7 &&  54 & 54 & 54 & 53         &&  76 & 76 & 71 & 65          &&  45 & 37 & 33 & 29 \\
     & 8 &&  50 & 51 & 55 & \oom       &&  71 & 71 & 67 & 65          &&  41 & 37 & 33 & 29 \\
    \midrule
     \multicolumn{17}{l}{MINRES, preconditioned with standard GS-MG} \\
    \midrule
     & 5 &&  64 & 167 & \dnc & \dnc    && 84 & 213 & \dnc & \dnc      &&  124 & 219 & \dnc & \dnc \\
    \bottomrule
    \end{tabular}
     \caption{ Iteration counts for the unit square}\label{tab:paramdom}
\end{table}

\begin{table}[p]
    \vspace{-1em}
    \centering
    \begin{tabular}{p{1.em}p{1.em}|p{.1cm}p{3.5em}p{2.5em}p{2.5em}|p{.1cm}p{3.5em}p{2.5em}p{2.5em}|p{.1cm}p{3.5em}p{2.5em}p{2.5em}}
    \toprule
          &          & \multicolumn{4}{c|}{Taylor-Hood} &  \multicolumn{4}{c|}{\Nedelec} & \multicolumn{4}{c}{Raviart-Thomas}\\
    $p$      &  $\ell$    &&  dof & $v$ & $p$ &&  dof & $v$ & $p$ &&  dof & $v$ & $p$   \\
    \midrule
    2    & 4 &&   2372        & 2e-5 & 1e-5 &&  1637     & 2e-5 & 4e-5 && 869        & 3e-4 & 3e-2 \\
          & 5 &&   9348        & 1e-6 & 6e-7 &&  6341      & 1e-6 & 4e-5 && 3269     & 3e-5 & 2e-2 \\
          & 6 &&   37124     & 7e-8 & 4e-6 &&  24965   & 7e-7 & 9e-5 && 12677    & 7e-6 & 2e-2 \\
          & 7 &&   147972   & 2e-8 & 7e-7 &&  99077   & 8e-7 & 1e-4 && 49925    & 3e-6 & 2e-2 \\
    \midrule
    5     & 4 &&    2891    & 2e-9 & 4e-8 &&  2066   & 6e-8 & 2e-6 && 1202     & 9e-7 & 6e-3 \\
          & 5 &&   10347    & 2e-9 & 1e-7 &&  7154   & 8e-8 & 5e-6 && 3890     & 1e-6 & 3e-3 \\
          & 6 &&   39083    & 3e-9 & 2e-7 &&  26546  & 9e-7 & 2e-4 && 13876   & 6e-7 & 3e-3\\
          & 7 &&   151951  & 7e-9 & 4e-7 &&  102194  & 2e-6 & 3e-4 && 52274 & 6e-7 & 4e-3\\
    \bottomrule
    \end{tabular}
     \caption{ Problem size and $L_2$-errors for the unit square}\label{tab:approx}
\end{table}

\begin{table}[p]
    \vspace{-1em}
    \centering 
    \begin{tabular}{p{.3em}p{2.2em}|p{.1cm}p{2em}p{2em}p{2em}p{2em}|p{.1cm}p{2em}p{2em}p{2em}p{2em}|p{.1cm}p{2em}p{2em}p{2em}p{2em}}
    \toprule
         &           & \multicolumn{5}{c|}{Taylor-Hood} &  \multicolumn{5}{c|}{\Nedelec} & \multicolumn{5}{c}{Raviart-Thomas} \\
    \quad & $\ell\,\diagdown\, p$     &&  2 & 3 & 5 & 8 && 2 & 3 & 5 & 8 && 2 & 3 & 5 & 8    \\
    \midrule
     \multicolumn{17}{l}{MINRES, preconditioned with SCMS-MG} \\
    \midrule
     & 5 &&  195 & 190 & 185 & 172         &&  257 & 246 & 244 & 206          &&  244 & 139 & 128 & 116 \\
     & 6 &&  208 & 217 & 213 & 199         &&  295 & 296 & 280 & 241          &&  192 & 170 & 142 & 129 \\
     & 7 &&  220 & 222 & 232 & 219         &&  329 & 330 & 314 & 281          &&  213 & 195 & 158 & 140 \\
     & 8 &&  231 & 239 & 244 & \oom       &&  333 & 342 & 333 & 306         &&  223 & 200 & 168 & 149 \\
    \midrule
     \multicolumn{17}{l}{MINRES, preconditioned with SCMS-MG-geo} \\
    \midrule
     & 5 &&  72 & 69 & 68 & 72             &&  69 & 69 & 65 & 63              &&  73 & 62 & 53 & 56 \\
     & 6 &&  77 & 75 & 73 & 79             &&  76 & 74 & 64 & 70              &&  71 & 69 & 59 & 63 \\
     & 7 &&  72 & 71 & 70 & 84             &&  79 & 70 & 68 & 74              &&  75 & 74 & 64 & 69 \\
     & 8 &&  74 & 73 & 72 & \oom           &&  73 & 73 & 71 & 78              &&  71 & 70 & 68 & 74 \\
    \midrule
     \multicolumn{17}{l}{MINRES, preconditioned with standard GS-MG} \\
    \midrule
     & 5 &&  70 & 173 & \dnc & \dnc        &&  110 & 225 & \dnc & \dnc       &&  182 & 220 & \dnc & \dnc \\
    \bottomrule
    \end{tabular}
     \caption{ Iteration counts for the quarter annulus (direct composition)}\label{tab:physdom}
\end{table}

\begin{table}[p]
    \vspace{-1em}
    \centering 
    \begin{tabular}{p{.3em}p{2.2em}|p{.1cm}p{2em}p{2em}p{2em}p{2em}|p{.1cm}p{2em}p{2em}p{2em}p{2em}|p{.1cm}p{2em}p{2em}p{2em}p{2em}}
    \toprule
          &          & \multicolumn{5}{c|}{Taylor-Hood} &  \multicolumn{5}{c|}{\Nedelec} & \multicolumn{5}{c}{Raviart-Thomas} \\
    \quad & $\ell\,\diagdown\, p$      &&  2 & 3 & 5 & 8 && 2 & 3 & 5 & 8 && 2 & 3 & 5 & 8    \\
    \midrule
     \multicolumn{17}{l}{MINRES, preconditioned with SCMS-MG} \\
    \midrule
     & 5 &&  331 & 331 & 338 & 317         &&  288 & 313 & 332 & 305          &&  480 & 309 & 295 & 300 \\
     & 6 &&  407 & 400 & 402 & 371         &&  361 & 387 & 405 & 374          &&  368 & 344 & 323 & 299 \\
     & 7 &&  452 & 455 & 455 & 450         &&  413 & 450 & 476 & 476          &&  418 & 395 & 367 & 341 \\
     & 8 &&  487 & 485 & 500 & \oom        &&  458 & 494 & 556 & 568          &&  441 & 438 & 411 & 361 \\
    \midrule
     \multicolumn{17}{l}{MINRES, preconditioned with standard GS-MG} \\
    \midrule
     & 5 &&  70 & 165 & \dnc & \dnc         &&  69 & 164 & \dnc & \dnc          &&  206 & 199 & \dnc & \dnc \\
    \bottomrule
    \end{tabular}
     \caption{ Iteration counts for the quarter annulus (Piola transform)}\label{tab:physdom:piola}
\end{table}

For the case of the quarter annulus, we distinguish between the results obtained
by the direct composition (Table~\ref{tab:physdom}) and for the Piola transform (Table~\ref{tab:physdom:piola}).
Again, we obtain first that GS-MG is robust in $h$, but the
convergence deteriorates if the polynomial degree grows. As it leads to better results, we have set up
the GS-MG on the physical domain. For the proposed SCMS-MG preconditioner,
observe that the results behave similar to the results for the unit square, however the iteration counts
are much larger, particularly if the Piola transform is used. For the direct composition, it is possible
to improve the convergence significantly by replacing the mass and stiffness
matrix on the parameter domain by a simple tensor-rank-one approximation of those matrices
on the physical domain (\emph{SCMS-MG-geo}). Note that the tensor-rank-one
approximation does not lead to any additional computational costs after the assembling phase.
The extension of such a rank-one geometry approximation to the
Piola transform is not yet known.  For the original SCMS-MG preconditioner,
have chosen $\beta$ and $\sigma$ as for the first model problem, just for the Raviart-Thomas smoother
for the case with Piola transformation, we have chosen $\beta=0.0025$.
For the rank-one corrected version, we have chosen $\beta = 0.01$; the damping has been chosen based on
an approximation for constants of the inverse inequality on the physical domain.

As in the case of standard finite elements, there are several possibilities to discretize the mixed formulation of the Stokes equations. Our
experiments indicate that it might pay off to use the (in terms of degrees of freedom) more expensive variant of Taylor Hood discretizations
than the other variants, particularly because it is known that this discretization also works for direct composition.
The $p$-robust smoothers which we have proposed for for the Poisson problem can be carried over also to the Stokes flow problem,
however it seems that further investigation is necessary concerning its application in the framework of non-trivial geometry transformations.

\bibliographystyle{abbrv}

\bibliography{paper}

\begin{thebibliography}{10}

\bibitem{BuffaEtAl}
G.~S. A.~Buffa, C.~de~Falco.
\newblock Isogeometric analysis: new stable elements for the stokes equation.
\newblock {\em Int.~J.~Num.~Meth.~Fluids}, 2010.

\bibitem{Bressan:2013}
A.~Bressan and G.~Sangalli.
\newblock Isogeometric discretizations of the {Stokes} problem: stability
  analysis by the macroelement technique.
\newblock {\em IMA Journal of Numerical Analysis}, 33(2):629--651, 2013.

\bibitem{Brezzi:1974}
F.~Brezzi.
\newblock {On the Existence, Uniqueness and Approximation of Saddle Point
  Problems Arising from Lagrangian Multipliers}.
\newblock {\em RAIRO Anal. Num\'{e}r.}, 8(2):129 -- 151, 1974.

\bibitem{Brezzi:Fortin:1991}
F.~Brezzi and M.~Fortin.
\newblock {\em {Mixed and Hybrid Finite Element Methods}}.
\newblock Springer-Verlag, 1991.

\bibitem{Buffa:2011}
A.~Buffa, C.~de~Falco, and G.~Sangalli.
\newblock Isogeometric analysis: Stable elements for the 2d {Stokes} equation.
\newblock {\em International Journal for Numerical Methods in Fluids},
  65(11-12):1407--1422, 2011.

\bibitem{DeBoor:1972}
C.~de~Boor.
\newblock On calculating with {B}-splines.
\newblock {\em Journal of Approximation Theory}, 6(1):50--62, 1972.

\bibitem{Evans:2013}
J.~A. Evans and T.~J.~R. Hughes.
\newblock Isogeometric divergence-conforming {B-splines} for the steady
  {Navier--Stokes} equations.
\newblock {\em Mathematical Models and Methods in Applied Sciences},
  23(08):1421--1478, 2013.

\bibitem{Greenbaum:1996}
A.~Greenbaum, V.~Pták, and Z.~Strakoš.
\newblock Any nonincreasing convergence curve is possible for {GMRES}.
\newblock {\em SIAM Journal on Matrix Analysis and Applications},
  17(3):465--469, 1996.

\bibitem{HT:2016}
C.~Hofreither and S.~Takacs.
\newblock Robust multigrid for isogeometric analysis based on stable splittings
  of spline spaces.
\newblock ArXiv e-print 1607.05035. \url{http://arxiv.org/abs/1607.05035}, July
  2016.

\bibitem{HTZ:2016}
C.~Hofreither, S.~Takacs, and W.~Zulehner.
\newblock A robust multigrid method for isogeometric analysis in two dimensions
  using boundary correction.
\newblock {\em Computer Methods in Applied Mechanics and Engineering}, 2016.
\newblock Available online.

\bibitem{Hughes:2005}
T.~J.~R. Hughes, J.~A. Cottrell, and Y.~Bazilevs.
\newblock Isogeometric analysis: {CAD}, finite elements, {NURBS}, exact
  geometry and mesh refinement.
\newblock {\em Computer Methods in Applied Mechanics and Engineering},
  194(39-41):4135--4195, Oct. 2005.

\bibitem{Larin:Reusken:2008}
M.~Larin and A.~Reusken.
\newblock A comparative study of efficient iterative solvers for generalized
  stokes problem.
\newblock {\em Numer. Linear Algebra Appl.}, 15:13--34, 2008.

\bibitem{gismoweb}
A.~Mantzaflaris, S.~Takacs, et~al.
\newblock {G+Smo (Geometry plus Simulation modules) v0.8.1}.
\newblock \url{http://gs.jku.at/gismo}, 2017.

\bibitem{Vanka:1986}
S.~P. Vanka.
\newblock {Block-implicit multigrid solution of Navier-Stokes equations in
  primitive variables}.
\newblock {\em Math. Comp.}, 65:138 -- 158, 1986.

\end{thebibliography}

\end{document}